\def \la {\lambda}

\def \tl {{\tilde {\lambda}}}
\def \tm {{\tilde {\mu}}}
\def \arr {\rightarrow}

\magnification=1200

\centerline {\bf A homological interpretation of Jantzen's sum formula}

\vskip 1cm

\centerline {\bf Upendra Kulkarni}

\vskip .3cm

\centerline {Truman State University, Kirksville, MO 63501  (Email: kulkarni@truman.edu)}

\vskip 1cm

\centerline{\bf Abstract}

\vskip .5cm

For a split reductive algebraic group, this paper observes a homological
interpretation for Weyl module multiplicities in Jantzen's sum formula.
This interpretation involves an Euler characteristic $\chi$ built from
$Ext$ groups between {\it integral} Weyl modules. The new interpretation makes
transparent For $GL_n$ (and conceivable for other classical groups) a certain
invariance of Jantzen's sum formula under ``Howe duality" in the sense of
Adamovich and Rybnikov. For $GL_n$ a simple and explicit general formula is
derived for $\chi$ between an arbitrary pair of integral Weyl modules. In light of
Brenti's work on certain $R$-polynomials, this formula raises interesting questions
about the possibility of relating $Ext$ groups between Weyl modules to
Kazhdan-Lusztig combinatorics.

\vskip 1cm

\noindent 0. {\it Introduction}

\vskip .5cm

Let $G_{\bf Z}$ be a split and connected reductive algebraic group
scheme over {\bf Z}. For a prime number $p$, $G_p$ will denote the
corresponding group scheme over ${\bf F}_p$, the field of $p$
elements. We will be concerned with (rational) representations
of $G_p$, in the course of which we will need to use $G_{\bf Z}$.
We will need some background material on these topics. For all
such standard facts see the recent edition of Jantzen's classic
text [Jantzen], where one can also find the original references.

\vskip .3cm

For a finite dimensional rational representation $M$ of $G_p$,
its formal character is
$$ch(M) = \sum_{\mu \in X} {\rm dim} (M_\mu) e(\mu) $$
in ${\bf Z}[X]^W$, the Weyl group invariants in the integral group
ring of the character group $X$ of a split maximal torus.
A central problem is to calculate formal characters of all
simple modules, which are in bijective correspondence with
their highest weights.
To get to the issue of interest in this paper, first let $L_p(\la)$
and $V_p(\la)$ respectively be the simple module and the Weyl module
(i.e., the universal highest weight module) corresponding to a
given dominant integral weight $\la$.
Since each of the families $\{ch(L_p(\la)\}$ and $\{ch(V_p(\la))\}$ forms
a basis of ${\bf Z}[X]^W$, and since $ch(V_p(\la))$ is known by Weyl's
character formula, the problem stated above reduces to finding the
multiplicities $a_{p,\la \mu}$ defined by
$$ch(L_p(\la)) = \sum_{{\rm dominant} \; \mu} a_{p, \la \mu} \; ch(V_p(\mu)).$$
There is the following well-known interpretaion of $a_{p,\la \mu}$ as an
Euler characteristic [Jantzen II.6.21]
$$a_{p, \la \mu} = \sum_i (-1)^i \; {\rm dim} \; {Ext^i}_{G_p} (V_p (\mu), L_p (\la)). \leqno (*)$$
(A lot more is known about this situation. Here is a sketch for completeness, though
we will not need this information later. To begin with, $a_{p, \la \la} = 1$.
Using the linkage principle, $a_{p, \la \mu} = 0$ unless $\mu \uparrow \la$. Using
translation functors further, we may additionally take without loss of generality
$\la$ and $\mu$ to be $p$-regular weights, provided the prime $p$ is large enough to
allow this. Then Lusztig's conjecture says that under a further condition on $\la$
(namely that it be in Jantzen's region) this multiplicity is the value at 1 of the
appropriate Kazhdan-Lusztig polynomial for the affine Weyl group  $W_p$ associated
to $G_p$. See [Jantzen II.6-8] for these issues and [Jantzen II.B,C,H] for
current status of Lusztig's conjecture.)

\vskip .3cm

The main purpose of this paper is twofold.
(1) First, to observe a homological interpretation in the spirit of (*) for
Weyl module multiplicities in Jantzen's sum formula [Jantzen II.8.19].
This formula calculates
$\sum_{i>0} ch(V^i_p(\la))$,
the sum of formal characters of all proper submodules $V^i_p(\la)$ appearing
in Jantzen's filtration of a Weyl module $V_p(\la)$. The new interpretation
involves an Euler characteristic $\chi$ over $G_{\bf Z}$ built from $Ext$ groups
between {\it integral} Weyl modules. (This is to be expected since Jantzen's
filtration is defined by working over {\bf Z}.) The proof rests on the same idea as
that behind (*), namely a fundamental $Ext$ calculation due to Cline-Parshall-Scott-van
der Kallen. Only here one uses the version over {\bf Z} rather than the one over
${\bf F}_p$.  An application is that for $GL_n$ the new interpretation
makes transparent a certain invariance of Jantzen's sum formula under ``Howe
duality" in the sense of Adamovich and Rybnikov. It further suggests that a similar
invariance may be true for some other classical groups too. The new interpretation
also provides some evidence for the likely importance in modular theory of $Ext$
groups over $G_{\bf Z}$, in particular those between integral Weyl modules. (The
modular analogue of $\chi$ is not so interesting for a pair of Weyl modules since
it is easily seen to satisfy the orthonormal property.)
With this in mind we will conduct additional analysis of the Euler characteristic
$\chi$ and briefly discuss some integral $Ext$ groups. All this is done in Section
1.

\vskip .3cm

(2) Second, to calculate for $GL_n$ a simple and explicit general formula for
the integral Euler characteristic $\chi$ between an arbitrary pair of integral
Weyl modules in terms of the associated dominant weights. This is done in Section
2. The formula is derived via a recursive procedure that employs the following
ingredients. Certain Weyl filtrations from [AB1] corresponding to characteristic-free
versions of Pieri-type rules, the skew representative theorem from [Kulkarni1],
conjugate symmetry of $Ext$ groups between Weyl modules from [AB2] and the
calculation in [Kulkarni2] of $Ext$ groups between Weyl modules for $GL_n$ whose
dominant weights differ by a single root.

\vskip .3cm

The $GL_n$ calculation of $\chi$ gives us $\sum_{i>0} ch(V^i_p(\la))$ in this
case, which is already known in general thanks to Jantzen's sum formula. Still
the $\chi$ calculation is of additional interest for the following reason.
The answer in Jantzen's formula is expressed as a linear combination of various
$ch(V_p(\mu))$, but due to the nature of this formula, it is not readily visible
which dominant weights $\mu$ occur. Moreover, the ones that do occur may repeat
(unless $\la$ is $p$-regular, see [Jantzen II.8.19, Remark 3]), so the multiplicities
of various $ch(V_p(\mu))$ cannot in general be read off easily from the formula.
Thus, especially for small primes, it is interesting to have a formula for $\chi$
between Weyl modules, which does give us these multiplicities.

\vskip 1cm

\noindent 1. {\it Jantzen's sum formula and $Ext$ groups}

\vskip .5cm

Let us proceed with the same set-up as in the introduction but the focus will
shift to representations of $G_{\bf Z}$ (= $G_{\bf Z}$-modules). Given a dominant
integral weight $\la$, both the Weyl module $V_p(\la)$ and its contravariant
dual $H_p(\la)$ (i.e., the dual Weyl module of largest weight $\la$) are
obtained from characteristic-free objects.
So we have {\bf Z}-free $G_{\bf Z}$-modules $V_{\bf Z}(\la)$ and $H_{\bf Z}(\la)$
such that $V_{\bf Z}(\la) \otimes {\bf F}_p = V_p(\la)$ and
$H_{\bf Z}(\la) \otimes {\bf F}_p = H_p(\la)$. One has
$ch (V_{\bf Z}(\la)) = ch (H_{\bf Z}(\la))  = ch (V_p(\la))$.
(The definition of $ch(M)$ in the introduction applies equally well to $G_{\bf Z}$-modules
that are free abelian groups of finite rank, the dimension being replaced by the
rank of each weight space as an abelian group.)

\vskip .3cm

$V_{\bf Z}(\la)$ and $H_{\bf Z}(\la)$ may be realized as follows. Extend the
scalars from {\bf Z} to rational numbers {\bf Q} to get the group scheme $G_{\bf Q}$.
Consider $V_{\bf Q}(\la)$, the simple $G_{\bf Q}$-module of highest weight $\la$.
Fix a vector $v$ in the one dimensional $\la$-weight space $V_{\bf Q}(\la)_{\la}$.
Among the finitely many $G_{\bf Z}$-stable lattices $M$ in $V_{\bf Q}(\la)$ such that
$M \cap V_{\bf Q}(\la)_{\la} = {\bf Z}v$, the unique minimal one is $V_{\bf Z}(\la)$
and the unique maximal one $H_{\bf Z}(\la)$. For future use note that this gives a
$G_{\bf Z}$-equivariant injection $\phi: V_{\bf Z}(\la) \hookrightarrow H_{\bf Z}(\la)$.
(In other words one has a nondegenerate bilinear form on $V_{\bf Z}(\la)$
(since as an abelian group $H_{\bf Z}(\la)$ is just the linear dual of $V_{\bf Z}(\la)$)
and this form is ``contravariant" due to the way $G_{\bf Z}$-action is defined on
$H_{\bf Z}(\la)$. See [Jantzen II.8.17]. We will not use this language here.)

\vskip .3cm

Prior to stating and proving the main results in 1.4, we will gather some general
results regarding $Ext$ groups (1.1), define and analyze the torsion Euler
characteristic $\chi$ (1.2) and review the setting of Jantzen's sum formula (1.3).
In view of the likely but as yet unclear significance of integral $Ext$ groups, some
additional remarks regarding $\chi$ and $Ext$ groups are offered in 1.5. Sections
1.6 discusses a connection of Jantzen's sum formula with Howe duality in the sense
of Adamovich and Rybnikov. Finally 1.7 points out possible generalizations to other
settings. All $G_{\bf Z}$-modules (except when scalars are extended to {\bf Q}) will
be finitely generated--equivalently, finitely generated as abelian groups.
For a finitely generated abelian group $M$, $M_{tor}$ will denote its torsion
subgroup and $M_{fr} = M/M_{tor}$ its largest {\bf Z}-free quotient. If $M$ is a
$G_{\bf Z}$-module so are $M_{tor}$ and $M_{fr}$.

\vskip .3cm

1.1. {\it Some basic results on $Ext$ groups over $G_{\bf Z}$} (two fundamental $Ext$
calculations, finiteness of $Ext$ groups in general).
Let us first record two important $Ext$ calculations due to Cline-Parshall-Scott-van
der Kallen.

\vskip .3cm

{\it Theorem.} [CPSvdK] For dominant weights $\la$ and $\mu$,

(i) $Ext^i_{G_{\bf Z}} (V_{\bf Z}(\mu), H_{\bf Z}(\la)) = 0$ unless ($\mu = \la$ and $i = 0$).

(ii) $Ext^i_{G_{\bf Z}} (V_{\bf Z}(\mu), V_{\bf Z}(\la)) = 0$ unless $\mu < \la$ or
($\mu = \la$ and $i = 0$).

{\it Note.} Of course $Hom_{G_{\bf Z}} (V_{\bf Z}(\la), V_{\bf Z}(\la)) \simeq
Hom_{G_{\bf Z}} (V_{\bf Z}(\la), H_{\bf Z}(\la))  \simeq {\bf Z}$.

\vskip .3cm

The next result proves finite generation of $Ext^i_{G_{\bf Z}}$ (and therefore
finiteness for $i > 1$). It is surely known (e.g., see [Jantzen II.4.10] for finite
dimensionality of $Ext^i_{G_p}$), but does not seem to be published, so a proof
is supplied. ({\it Note.} A more basic and immediate argument for finite generation
than the one given is as follows. Use an appropriate Schur algebra $S$ over {\bf Z},
which is finitely generated, and appeal to Donkin's theorem (actually to its easy
extension to the integral situation) asserting isomorphism of $Ext_S$ with
$Ext_{G_{\bf Z}}$. Instead a different argument is given below to illustrate two
techniques: use of Theorem 1.1(i) and a certain strategy to build up general
$G_{\bf Z}$-modules.)

\vskip .3cm

{\it Proposition.} Fix finitely generated $G_{\bf Z}$-modules $M$ and $N$ and consider the
abelian groups $Ext^i_{G_{\bf Z}}(M,N)$.

(i)   These groups are all finitely generated and are $0$ for large $i$.

(ii)  For $i>0$ these groups are all finite.

(iii) $Hom_{G_{\bf Z}}(M,N)$ is infinite iff the $G_{\bf Q}$-modules $M  \otimes {\bf Q}$
and $N \otimes {\bf Q}$ have at least one isomorphic simple summand.

\vskip .3cm

{\it Proof.}
(i) If the statement is true for $M$ = two of the three modules in a short exact sequence
then it is clearly true for the third. (We will say that each module is ``obtainable"
from the other two.) The same is true for $N$. We will show that all finitely generated
modules are obtainable from the family  $\{V_{\bf Z}(\la)\}$ as well as from $\{H_{\bf Z}(\la)\}$.
Now taking $M = V_{\bf Z}(\mu)$ and $N = H_{\bf Z}(\la)$ will give (i) by Theorem 1.1(i)
quoted above.

\vskip .3cm

$M$ is obtainable from $M_{tor}$ and $M_{fr}$. Now $M_{tor}$ has a composition series
with factors $L_p(\la)$ for various primes $p$ and dominant weights $\la$. For fixed
$p$ the family $\{L_p(\la)\}$ is well-known to be obtainable from either of the families
$\{V_p(\la)\}$ and $\{H_p(\la)\}$, which in turn are respectively obtainable from
the families $\{V_{\bf Z}(\la)\}$ and $\{H_{\bf Z}(\la)\}$ using multiplication by
$p$. As for $M_{fr}$, it has a filtration whose factors are $G_{\bf Z}$-stable lattices
$P$ in various $V_{\bf Q}(\la)$. We have lattices of the same rank
$V_{\bf Z}(\la) \subset P \subset H_{\bf Z}(\la)$, resulting in torsion modules
$P/V_{\bf Z}(\la)$ and $H_{\bf Z}(\la)/P$. By the argument for torsion modules $M_{fr}$
is also obtainable from either of the families $\{V_{\bf Z}(\la)\}$ and $\{H_{\bf Z}(\la)\}$.
This proves (i).

\vskip .3cm

(ii) In view of (i) clearly all $Ext^i$ are finite if $M$ or $N$ is torsion. So by looking
at the appropriate long exact sequences derived from $0 \arr M_{tor} \arr M \arr M_{fr}$
and $0 \arr N_{tor} \arr N \arr N_{fr}$ it suffices to take both $M$ and $N$ to be {\bf Z}-free.
Now the result follows after extending scalars to {\bf Q} (e.g., use the universal
coefficient theorem [Jantzen I.4.18] to relate $Ext_{G_{\bf Z}}$ and $Ext_{G_{\bf Q}}$), where
all representations become semisimple.

\vskip .3cm

Finally (iii) is clear.

\vskip .3cm

1.2. {\it The torsion Euler characteristic $\chi$.} In view of Proposition 1.1 the following
definition is valid.

\vskip .3cm

{\it Definition.} For finitely generated $G_{\bf Z}$-modules $M$ and $N$ define
$$\chi(M,N) = |Hom_{G_{\bf Z}}(M,N)_{tor}| \prod_{i > 0} |Ext^i_{G_{\bf Z}}(M,N)|^{(-1)^i}.$$
{\it Note.} Another alternative is to simply drop the first term if the $Hom$ is
infinite. In fact for the purposes of this paper that would work just as well, since
in all our uses of $\chi$, $Hom_{G_{\bf Z}}(M,N)$ will be either {\bf Z}-free or
torsion. But in the general definition it seems preferable to include the ``finite
part" $Hom_{G_{\bf Z}}(M,N)_{tor}$.

\vskip .3cm

Clearly $\chi$ is multiplicative on short exact sequences of finite modules in either
argument. Moreover if one of the arguments of $\chi$ is finite, $\chi$ is multiplicative
on short exact sequences in the other argument. This implies that for torsion modules
$M$ and $N$, $\chi(M,N)$ is always 1. ({\it Proof.} Apply $Hom_{G_{\bf Z}}(M,-)$ to
the exact sequence
$0 \arr V_{\bf Z}(\nu) \buildrel p \over \arr V_{\bf Z}(\nu) \arr V_p(\nu) \arr 0$
to get $\chi(M, V_p(\nu))=1$. Now any finite $N$ may be obtained from the various
$V_p(\nu)$.)

\vskip .3cm

$\chi$ fails to be multiplicative in general because infinite $Hom$ groups exist. Let us sketch
how this failure can be quantified. Apply $Hom_{G_{\bf Z}}(M, -)$ or $Hom_{G_{\bf Z}}(-, M)$
to a short exact sequence $0 \arr P \arr Q \arr R \arr 0$ of $G_{\bf Z}$-modules and take
the corresponding long exact sequence. Replacing the three $Hom$ terms by their torsion
parts renders the long sequence possibly inexact in two places: the last $Hom$ term and the
first $Ext^1$ term. Let us make explicit the necessary adjustment to multiplicativity in two
situations of interest.

\vskip .3cm

{\it Example 1.} (Breaking up $M,N$ into free and torsion parts.) Apply
$Hom_{G_{\bf Z}}(-,N_{tor})$ to $0 \arr M_{tor}\arr M \arr M_{fr}$ to get
$\chi(M, N_{tor}) = \chi(M_{fr}, N_{tor})$ since $\chi(M_{tor}, N_{tor}) = 1$. Apply
$Hom_{G_{\bf Z}} (-,N_{fr})$ to the same sequence and use $Hom_{G_{\bf Z}}(M_{tor}, N_{fr})=0$
to get $\chi(M, N_{fr}) = \chi(M_{tor}, N_{fr}) \chi(M_{fr}, N_{fr})$.
Next apply $Hom_{G_{\bf Z}}(M,-)$ to  $0 \arr N_{tor}\arr N \arr N_{fr} \arr 0$ to get
a long exact sequence beginning as follows.
$$0 \arr Hom_{G_{\bf Z}}(M,N_{tor}) \arr  Hom_{G_{\bf Z}}(M,N) \arr Hom_{G_{\bf Z}}(M,N_{fr})
\arr  Ext^1_{G_{\bf Z}}(M,N_{tor}) \arr \cdots $$
Now $Hom_{G_{\bf Z}}(M,N_{tor}) \simeq Hom_{G_{\bf Z}}(M,N)_{tor}$, giving an injection of free
abelian groups of equal rank $Hom_{G_{\bf Z}}(M,N)_{fr} \hookrightarrow Hom_{G_{\bf Z}}(M,N_{fr})$.
Calling the necessarily finite cardinality of the cokernel of this injection $s$,
we clearly have $\chi(M,N) = \chi(M,N_{tor}) \chi(M,N_{fr}) s $. All in all we have
$$\chi(M,N) = \chi(M_{fr}, N_{tor}) \chi(M_{tor}, N_{fr}) \chi(M_{fr}, N_{fr}) s.$$
It is easy to see that the correction factor $s$ is really necessary. Otherwise the map
$Hom_{G_{\bf Z}}(M,N) \arr Hom_{G_{\bf Z}}(M,N_{fr})$ would always be surjective,
i.e., the exact sequence $0 \arr N_{tor}\arr N \arr N_{fr} \arr 0$ would always split (e.g.,
letting $M=N_{fr}$, the inverse image of $id_{N_{fr}}$ would give a splitting), which surely
does not happen. To give a concrete example, it suffices to produce a non-split extension
of a {\bf Z}-free $G_{\bf Z}$-module by a torsion $G_{\bf Z}$-module. For this take
$G_{\bf Z} = GL({\bf Z}^n)$. Let $\Lambda^2$ and $D^2$ respectively be the second
exterior and divided powers of the defining representation ${\bf Z}^n$. One knows from
[AB2, Section 9] that $Ext^i_{G_{\bf Z}}(\Lambda^2, D^2)$ is ${\bf Z}/2{\bf Z}$ if $i = 1$
and vanishes for other $i$. Using this in the long exact sequence obtained by applying
$Hom_{G_{\bf Z}}(\Lambda^2, -)$
to $0 \arr D^2 \buildrel 2 \over \arr D^2 \arr D^2 \otimes {\bf Z}/2{\bf Z} \arr 0$
yields that $Ext^i_{G_{\bf Z}}(\Lambda^2, D^2 \otimes {\bf Z}/2{\bf Z}) = {\bf Z}/2{\bf Z}$
if $i = 0,1$ and vanishes otherwise. In particular we get the desired non-split extension.
See Remark 1.5.(5) below for similar elementary analysis of the relationships among
several kinds of $Ext$ groups in a general situation.

\vskip .3cm

{\it Example 2.} ($\chi$ involving an extension of {\bf Z}-free $G_{\bf Z}$-modules.)
Take an exact sequence of {\bf Z}-free $G_{\bf Z}$-modules $0 \arr N \arr P \arr R \arr 0.$
Then for a {\bf Z}-free $G_{\bf Z}$-module $M$ one has exactly as in Example 1
$$\chi(M,P) = \chi(M,N) \chi(M,R)  s ,$$
where $s$ is the cardinality of the necessarily finite cokernel of the map
of free abelian groups $Hom (M,P) \arr Hom(M,R)$. Clearly $s=1$ unless
$M \otimes {\bf Q}, P \otimes {\bf Q},  R \otimes {\bf Q}$ all have an
isomorphic simple summand. Similar considerations apply to such a short exact sequence
in the first argument of $\chi$. In the recursive algorithm to compute $\chi$ between
integral Weyl modules for $GL_n$ to be presented in the next section, the terminating
stage of the algorithm involves exactly such pairs $(M,R)$ which have infinite $Hom$
groups. In this sense these numbers $s$ are the source of the final numerical answers
for $\chi$.

\vskip .3cm

1.3. {\it Background on formal character for torsion modules and Jantzen's sum formula.}
We need some more machinery to be able to state and prove the main result. A convenient
reference for all of this background material is [Jantzen II.8].

\vskip .3cm

{\it Definition and basic properties of $ch_{tor}$.}
One can define a torsion formal character for a finite $G_{\bf Z}$-module $M$ by
$$ch_{tor}(M) = \sum_{\mu \in X} div|M_\mu| e(\mu),$$
where $div$ stands for taking the divisor of a rational number (here an integer).
(Compare the definition of $\nu^c$ in [Jantzen II.8.12].) Clearly
$ch_{tor}(L_p(\la)) = ch (L_p(\la)) [p]$ and $ch_{tor}$ is additive on short
exact sequences of finite $G_{\bf Z}$-modules. As $\la$ varies over all dominant
weights and $p$ over all primes, each of the two families
$\{ch_{tor}(L_p(\la))\}$ and $\{ch_{tor}(V_p(\la)) = ch_{tor} (H_p(\la)) \}$
forms a basis of the abelian group of all torsion formal characters. Note that since
$ch(V_p(\la)) = ch(V_{\bf Z}(\la))$ is independent of $p$, it makes sense to speak
about the coefficient of $ch(V_{\bf Z}(\la))$ in $ch_{tor}(M)$, this coefficient
being the divisor of a unique positive rational number (that will be 1 for all but
finitely many $\la$). We will need these considerations while explaining the
setting of Jantzen's sum formula, which will be our next task.

\vskip .3cm

{\it Jantzen's filtration and Jantzen's sum formula.}
Recall the injection $\phi: V_{\bf Z}(\la) \hookrightarrow H_{\bf Z}(\la)$.
Jatzen's filtration is a descending filtration $V^i_p(\la)$ of $V_p(\la)$ defined
as follows. Fixing $p$ for the moment, first let $V^i_{\bf Z}(\la)$ be the submodule
$\phi^{-1} (p^i H_{\bf Z}(\la))$ of $V_{\bf Z}(\la)$ and then $V^i_p(\la) =$ the
image of $V^i_{\bf Z}(\la)$ under the canonical map $V_{\bf Z}(\la) \arr V_p(\la)$.
(So in particular $V_p(\la)/V^1_p(\la)$ = the image of $\phi \otimes id_{{\bf F}_p}$,
which is well-known to be $L_p(\la)$.)

\vskip .3cm

Set $Q(\la) = coker(\phi)$. By a well-known argument due to Jantzen,
$$ch_{tor} Q(\la) =  \sum_p \left ( \sum_{i>0} ch(V^i_p(\la)) \right ) [p].$$
At the same time it is immediate (e.g., after diagonalizing $\phi$ by
using suitable bases for $V_{\bf Z}(\la)$ and $H_{\bf Z}(\la)$) that calculating
$ch_{tor} Q(\la)$ is equivalent to calculating the determinant of $\phi$ on each
weight space. In most cases Jantzen succeeded in calculating the determinant
(more precisely its $p$-adic valuation when $p$ is not small and without restriction
for type A) and so proved a formula for $\sum_{i>0} ch(V^i_p(\la))$, which was
soon afterwards obtained in general by Andersen via a different method. See
[Jantzen II.8] for a detailed explanation of all this and the formula itself.
We do not need the formula here. Rather, our goal is to observe an Euler
characteristic interpretation for the Weyl module multiplicities occuring in
this formula. These multiplicities, denoted $b_{p, \la \mu}$, are defined by
$$\sum_{i>0} ch(V^i_p(\la)) = \sum_{{\rm dominant} \, \mu} b_{p,\la \mu} \; ch(V_{\bf Z}(\mu)).$$
Note that only $\mu < \la$ can have nonzero coefficients on the right hand side.

\vskip .3cm

1.4. {\it A formula for $ch_{tor}$ using $\chi$ and application to the sum formula.}
We will interpret the coefficients $b_{p, \la \mu}$ in terms of $\chi$ between
integral Weyl modules. This will be a consequence of the following integral
analogue of (*) in the introduction (more precisely, analogue of [Jantzen II.6.21]).

\vskip .3cm

{\it Proposition.} For a finite $G_{\bf Z}$-module $M$,
$$\eqalign{
ch_{tor}(M) & \quad =  \sum_{{\rm dominant} \, \mu}   div(\chi(V_{\bf Z}(\mu),M)) \; ch(V_{\bf Z}(\mu)) \cr
            & \quad = \sum_{{\rm dominant} \, \mu} - div(\chi(M,H_{\bf Z}(\mu))) \; ch(V_{\bf Z}(\mu)).\cr
}$$
Further, the statements stay true if one replaces $V_{\bf Z}(\mu)$ or $H_{\bf Z}(\mu)$
by any $G_{\bf Z}$-stable lattice inside $V_{\bf Q}(\mu)$.

\vskip .3cm

{\it Proof.} By additivity of $ch_{tor}$ and of the divisor of $\chi(V_{\bf Z}(\mu),-)$
on finite $G_{\bf Z}$-modules it is enough to check the first equality for $M = H_p(\nu)$
for each prime $p$ and each dominant integral weight $\nu$. One has
$ch_{tor}(H_p(\nu)) = ch(V_{\bf Z}(\nu)) [p]$. For the other calculation apply
$Hom_{G_{\bf Z}}(V_{\bf Z}(\mu), -)$ to
$0 \arr H_{\bf Z}(\nu) \buildrel p \over \arr H_{\bf Z}(\nu) \arr H_p(\nu) \arr 0$
and use Theorem 1.1(i) above to get $\chi(V_{\bf Z}(\mu), H_p(\nu)) = p ^{\delta_{\mu \nu}}$.
The second equality is checked similarly using
$M = V_p(\nu)$. Since $\chi$ is 1 for a pair of finite $G_{\bf Z}$-modules, using
multiplicativity of $\chi(-,M)$ (respectively, of $\chi(M,-)$) one may replace $V_{\bf Z}(\mu)$
(respectively $H_{\bf Z}(\mu)$) inside $\chi$ by any $G_{\bf Z}$-stable lattice inside
$V_{\bf Q}(\mu)$, thus proving the last statement.

\vskip .3cm

{\it Corollary.} (Homological interpretation of Jantzen's sum formula.) The
multiplicities of $ch(V_{\bf Z}(\mu))$ in Jantzen's sum formulas for $V_p(\la)$
for various primes $p$ are given by
$$\sum_p b_{p, \la \mu} [p] = - div (\chi(V_{\bf Z}(\mu), V_{\bf Z}(\la)).$$

{\it Proof.} Note that
$$\eqalign {
ch_{tor} Q(\la)  & = \sum_p  \left (\sum_{{\rm dominant} \, \mu}
                                       b_{p, \la \mu} ch(V_{\bf Z}(\mu)) \right) [p]
                   = \sum_{{\rm dominant} \, \mu}
                             \left (\sum_p b_{p, \la \mu} [p] \right) ch(V_{\bf Z}(\mu)) \cr
                 &  = \sum_{{\rm dominant} \, \mu} div(\chi(V_{\bf Z}(\mu),Q(\la))) \; ch(V_{\bf Z}(\mu)).\cr
}$$
Here the first line comes simply from the set-up in 1.3 and the second from the first
equality in Proposition 1.4. So all we need to show is that
$$div (\chi(V_{\bf Z}(\mu), Q(\la))) = - div (\chi(V_{\bf Z}(\mu), V_{\bf Z}(\la))).$$
For this take the long exact sequence obtained by applying $Hom_{G_{\bf Z}}(V_{\bf Z}(\mu), -)$
to $0 \arr V_{\bf Z}(\la) \arr H_{\bf Z}(\la) \arr Q(\la) \arr 0$. Now
$\chi(V_{\bf Z}(\mu), H_{\bf Z}(\la)) = 1$ by using Theorem 1.1(i) once again and hence
we have the desired result. (Note that the first two $Hom$ terms in the long exact
sequence--potentially the only infinite ones--vanish unless $\mu = \la$. If $\mu = \la$,
these terms are {\bf Z} and the map between them is an isomorphism. So multiplicativity
of $\chi(V_{\bf Z}(\mu), -)$ is valid in either case. The latter case is trivial anyway,
since then all terms other than the first two vanish and all three $\chi$'s are 1.)

\vskip .3cm

1.5. {\it Complementary remarks on $\chi$ and on $Ext$ groups.}
(1) Note that the analogue over ${\bf F}_p$ of $\chi$ between Weyl modules is uninteresting as it
satisfies the orthonormal property. One has
$\sum_i (-1)^i {\rm dim} \; Ext^i_{G_p} (V_p(\la), V_p(\mu)) = \delta_{\la \mu}$
by [Jantzen II.6.21].

\vskip .3cm

(2) {\it $\chi$ for some other pairs of modules.}
One may contemplate $\chi(A,B)$ where $A$ and $B$ are $G_{\bf Z}$-modules from the
three families $\{V_{\bf Z}(\la)\}, \{H_{\bf Z}(\la)\}$ and $\{Q(\la)\}$. Using information
obtained so far, it is straightforward to analyze the nine possibilities.
By Theorem 1.1(i), $\chi(V_{\bf Z}(\mu), H_{\bf Z}(\la)) = 1$.
$\chi(Q(\mu), Q(\la)) = 1$ since both arguments are torsion modules.
$\chi(V_{\bf Z}(\mu), V_{\bf Z}(\la))
= \chi(H_{\bf Z}(\la), H_{\bf Z}(\mu))$ (by contravariant duality)
$= 1 / \chi(V_{\bf Z}(\mu), Q(\la))$ (by Theorem 1.1(i))
$= \chi(Q(\la), H_{\bf Z}(\mu))
= 1 / \chi(H_{\bf Z}(\mu), Q(\la))
= \chi(Q(\la), V_{\bf Z}(\mu))$ (last three expressions by Proposition 1.4).
This leaves $\chi(H_{\bf Z}(\mu), V_{\bf Z}(\la))$. To relate this to previous cases apply
$Hom_{G_{\bf Z}}(-, V_{\bf Z}(\la))$
to  $0 \arr V_{\bf Z}(\mu) \arr H_{\bf Z}(\mu) \arr Q(\mu) \arr 0$ and use known information
to get $\chi(H_{\bf Z}(\mu), V_{\bf Z}(\la)) =
\chi(V_{\bf Z}(\mu), V_{\bf Z}(\la)) \chi(V_{\bf Z}(\la), V_{\bf Z}(\mu))$. The symmetry in
$\la$ and $\mu$ is to be expected in view of contravariant duality. Note that by Theorem 1.1(ii)
at most one of the two factors on the right hand side may be different from 1.

\vskip .3cm

(3) {\it A certain skew-symmetry of $\chi$}. It follows from Proposition 1.4 that
for a torsion $G_{\bf Z}$-module $M$, $\chi(M,V_{\bf Z}(\mu)) \chi(V_{\bf Z}(\mu),M) = 1$.
Using multiplicativity of $\chi(M,-)$ as well as of $\chi(-,M)$, and since $\chi$ is 1
for a pair of torsion modules, one may replace $V_{\bf Z}(\mu)$ by any $G_{\bf Z}$-module
$N$ via a route similar to the one followed in the proof of Proposition 1.1(i). So for a
torsion $G_{\bf Z}$-module $M$ and for any $G_{\bf Z}$-module $N$, one has
$\chi(M,N) \chi(N,M) = 1$. This is not true in general if $M$ is not torsion, e.g., when
$M$ and $N$ are both integral Weyl modules in view of Corollary 1.4 and Theorem 1.1(ii).

\vskip .3cm

(4) Of course it is much more interesting (and harder) to calculate the $Ext$ groups
themselves rather than merely calculating $\chi$. Here is a reason why the $Ext$ groups
between integral Weyl modules are likely to be important. Note that we made essential
use of Theorem 1.1(i) only ``up to Euler characteristic," not in its full strength.
This theorem gives $Ext^{i+1}_{G_{\bf Z}}(V_{\bf Z}(\mu),V_{\bf Z}(\la)) \simeq
Ext^{i}_{G_{\bf Z}}(V_{\bf Z}(\mu),Q(\la))$ and clearly the structure of $Q(\la)$ is
intimately related to Jantzen's filtration itself. It would be very interesting if one
can relate individual $Ext$ groups between integral Weyl modules to Jantzen's filtration.
It seems reasonable to speculate that Kazhdan-Lusztig type combinatorics should come
into play in such a relationship. (There is already a hint to this effect for type A.
See Remark 2 after the proof of Theorem 2.3.) For instance one can ask the following
question. Consider the alternating sum obtained by taking the divisor of the defining
expression for $\chi(V_{\bf Z}(\mu),V_{\bf Z}(\la))$. What is the polynomial obtained
by replacing $-1$ in this expression by an indeterminate $q$? (Unlike over ${\bf F}_p$,
this is a weaker question than knowing the individual $Ext$ groups since in general
the structure of an abelian group is not determined by its size.)

\vskip .3cm

(5) {\it Elementary observations relating $Ext_{G_{\bf Z}}$ and $Ext_{G_p}$.}
Consider a $G_p$-module $L$ and a {\bf Z}-free ${G_{\bf Z}}$-module $M$. First we have
the isomorphism of ${\bf F}_p$-vector spaces (compare [McNinch, Lemma 3.1.1b])
$Ext^i_{G_{\bf Z}}(M,L) \simeq Ext^i_{G_p}(M \otimes {\bf F}_p,L)$.
{\it Proof.} Using $^*$ to denote the linear dual over {\bf Z}, [Jantzen I.4.2(1) and
I.4.4] give $Ext^i_{G_{\bf Z}}(M,L) \simeq H^i({G_{\bf Z}},M^* \otimes L)$ and
$Ext^i_{G_p}(M \otimes {\bf F}_p,L) \simeq H^i({G_p},M^* \otimes {\bf F}_p \otimes L)$.
Now the Hochschild complexes computing both group cohomologies are isomorphic.
Alternatively, take a projective resolution $P. \arr M$ over an appropriate Schur algebra
$S_{\bf Z}$. Tensoring by ${\bf F}_p$ gives a projective resolution of $M \otimes {\bf F}_p$
over the Schur algebra $S_p = S_{\bf Z} \otimes {\bf F}_p$. Applying the appropriate
$Hom$ to each resolution gives isomorphic complexes computing the required $Ext$ groups.

\vskip .3cm

Next, it is easy to relate  these isomorphic $Ext$ groups to
$Ext^i_{G_{\bf Z}}(M \otimes {\bf F}_p,L)$. Apply $Hom_{G_{\bf Z}}(-,L)$ to the
short exact sequence $0 \arr M \buildrel p \over \arr M \arr M \otimes {\bf F}_p \arr 0$.
The resulting long exact sequence easily breaks up into short exact sequences, giving
the isomorphism of ${\bf F}_p$-vector spaces $Ext^i_{G_{\bf Z}}(M \otimes {\bf F}_p,L)
                            \simeq Ext^{i-1}_{G_{\bf Z}}(M,L) \oplus Ext^i_{G_{\bf Z}}(M,L)$.

\vskip .3cm

One can say more in the following special situation. Let $L = N \otimes {\bf F}_p$, where $N$
is a {\bf Z}-free ${G_{\bf Z}}$-module. Suppose the localization at prime $p$ of the abelian
group $Ext^i_{G_{\bf Z}}(M,N)$ is a direct sum of $n_i$ cyclic groups. Using the
short exact sequences $0 \arr M \buildrel p \over \arr M \arr M \otimes {\bf F}_p \arr 0$
and $0 \arr N \buildrel p \over \arr N \arr N \otimes {\bf F}_p \arr 0$ as in the previous
paragraph, one gets the following equalities.
dim $Ext^i_{G_{\bf Z}}(M \otimes {\bf F}_p,N) = n_{i-1}+n_i$.
dim $Ext^i_{G_{\bf Z}}(M, N \otimes {\bf F}_p) =$
dim $Ext^i_{G_p}(M \otimes {\bf F}_p,N \otimes {\bf F}_p) = n_i + n_{i+1}$.
(The preceding equality can also be seen from the Universal Coefficient Theorem [Jantzen I.4.18a].)
dim $Ext^i_{G_{\bf Z}}(M \otimes {\bf F}_p,N \otimes {\bf F}_p) = n_{i-1}+2n_i+n_{i+1}$.
\vskip .3cm

The foregoing considerations may be of particular interest in two cases because of the
significance of the $Ext$ groups involved. Take $M = V_{\bf Z}(\mu)$ in each case
and in turn let $L = L_p(\la)$ (see (*) in the Introduction) or $L = V_p(\la)$ (in view
of the previous remark).

\vskip .3cm

1.6. {\it Symmetry of sum formula under Howe duality and complements.}
Consider Young diagrams of two partitions $\la$ and $\mu$. These may be considered
to be two dominant weights for $GL(F)$, where $F$ is a free abelian group of rank at
least as much as the number of rows in $\la$ as well as $\mu$. Let $\tl$ (respectively,
$\tm$) be the Young diagram obtained by transposing the rows and columns of $\la$
(respectively, $\mu$). $\tl,\tm$ are dominant weights for $GL(E)$, where $E$ is a free
abelian group of rank at least as much as the number of rows in $\tl$ as well as $\tm$.
(One could also work with a single group of large enough rank, or in the category of
polynomial functors without having to worry about the rank at all.) Now the functor
$\Omega$ from [AB2, Section 7] (a characteristic-free form of Howe duality,
in the sense similar to [AR]) combined with contravariant duality gives
the following ``conjugate symmetry" of $Ext$ groups between Weyl modules.
$$Ext^i_{GL(F)}(V_{\bf Z}(\mu),V_{\bf Z}(\la)) \simeq
                          Ext^i_{GL(E)}(V_{\bf Z}(\tl),V_{\bf Z}(\tm)).\leqno(1.6.1)$$
Still working with $GL(F)$, a more elementary symmetry of these $Ext$ groups is as
follows. Enclose Young diagrams of $\la$ and $\mu$ in a rectangle of height $rank(F)$
and of suitable width $w$. Let $\la^c$ and $\mu^c$ respectively be the partitions
obtained by taking the complements of $\la$ and $\mu$ within this rectangle. One has
(e.g., by [ABW, II.4]) $V_{\bf Z}(\mu) = H_{\bf Z}(\mu^c)^* \otimes (det)^w$ and
likewise for $\la$. (Here $^*$  denotes ordinary linear dual of a representation, where
one uses the group antiautomorphism taking $g$ to $g^{-1}$ to convert the natural
right action into a left one.) Now
$$\eqalign {Ext^i_{GL(F)}(V_{\bf Z}(\mu),V_{\bf Z}(\la))
&  \simeq Ext^i_{GL(F)}(H_{\bf Z}(\mu^c)^* \otimes (det)^w, H_{\bf Z}(\la^c)^* \otimes (det)^w) \cr
&  \simeq Ext^i_{GL(F)}(H_{\bf Z}(\mu^c)^*, H_{\bf Z}(\la^c)^*) \cr
&  \simeq Ext^i_{GL(F)}(H_{\bf Z}(\la^c), H_{\bf Z}(\mu^c)) \cr
&  \simeq Ext^i_{GL(F)}(V_{\bf Z}(\mu^c), V_{\bf Z}(\la^c)). \cr} \leqno (1.6.2)$$
Here the second line comes from canceling the determinant, the third by using
linear duality and the last by using contravariant duality.

\vskip .3cm

It follows immediately from (1.6.1), (1.6.2) and Corollary 1.4 that for general linear
groups Jantzen's sum formula is stable under Howe duality and under complements.
More precisely one has the following.

\vskip .3cm

{\it Corollary}. For the general linear groups (in the terminology of 1.3),
$$b_{p, \la \mu} = b_{p, \tm \tl} = b_{p, \la^c \mu^c}.$$

Adamovich and Rybnikov [AR] have constructed Howe duality functors in positive
characteristic for several other pairs of classical groups. Their set-up gives
an isomorphim of $Ext$ groups analogous to a combination of (1.6.1) and (1.6.2),
but in characteristic $p$. It is not clear that their isomorphism is valid while
working over {\bf Z}. The groups $Ext^i_{G_{\bf Z}}(V_{\bf Z}(\mu),V_{\bf Z}(\la))$
for all $i$ determine the corresponding modular $Ext$ groups (via, e.g.,
[Jantzen I.4.18a]), but not vice versa. Still the work of Adamovich and Rybnikov
makes it reasonable to ask whether Jantzen's sum formulas are stable under Howe
duality for any other pairs of groups besides the general linear groups. This should
(at least in principle) be verifiable directly from the known sum formulas.

\vskip .3cm

{\it Notes.} (1) [McNinch] proves the following intersting connection within this cluster
of ideas. Howe duality in [AR] carries Jantzen's filtrations to Andersen's ``tilting
filtrations."  (2) Here is yet another connection between Andersen's tilting filtartions
and Jantzen's filtrations. (This is entirely independent of Howe duality, but is still
mentioned here for completeness.) Using the results in this section, it is possible to
derive Andersen's ``titling sum formula" as a formal consequence of Jantzen's sum formula.
Currently this is known only if the characteristic is not small, since the original proof
uses regular weights. The new proof will be explained in [Kulkarni3].

\vskip .3cm

1.7. {\it Generalization to highest weight categories.}
The arguments in this section used certain objects (such as simple, Weyl and dual Weyl
modules indexed by a suitable partially ordered set) and facts about these objects
(such as suitable vanishing properties) that are available in other situations as well.
The natural home for the needed objects and vanishing results is the axiomatic notion
of a highest weight category (= representations of a quasi-hereditary algebra)
due to Cline-Parshall-Scott. Additionally one needs a suitable model of such a category
over a principal ideal domain from which the category of interest is obtained by
reduction modulo a prime. See for example the set-up in [McNinch]. It should be
straighforward to carry over the results in this section (except 1.6) to such an
axiomatic setting. In particular analogues of the main results should hold in the
following situations: the BGG cateory $\cal \char'117$ (using the version constructed by
Gabber-Joseph) and representations of quantum groups at a root of unity. (One has
Jantzen's filtrations and sum formulas in these situations as well.) Because of the
focused nature of interest in this paper, we will not carry out any of these
extensions here.

\vskip 1cm

\noindent 2. {\it A formula for $\chi(V_{\bf Z}(\mu), V_{\bf Z}(\la))$ for the general linear group}

\vskip .5cm

Henceforth let $G_{\bf Z} = GL(F)$, where $F$ is a free abelian group of finite rank.
The rank of $F$ will be essentially immaterial (see below for the precise statement).
We will calculate $\chi(V_{\bf Z}(\mu), V_{\bf Z}(\la))$ for arbitrary dominant weights
$\la$ and $\mu$ and so obtain for this group a version of Jantzen's sum formula (in
view of Corollary 1.4). Since the argument has some combinatorial intricacy, it will
be broken up as follows. In 2.1 we will develop an algorithm to calculate the desired
$\chi$. The algorithm will be illustrated in 2.2 by working out two simple classes of
examples, after which we will derive the general formula in 2.3. But first let us
(slightly) reformulate the problem so as to arrive at the setting in which the
algorithm will take place.

\vskip .3cm

Each of the dominant weights $\la$ and $\mu$ may be identified with a weakly decreasing
set of $rank(F)$ integers. By tensoring with the determinant enough times, we may take
these integers to be nonnegative without loss of generality. Thus $\la$ and $\mu$ may
be taken to be just partitions with at most $rank(F)$ rows. We will freely identify
partitions with their Young diagrams.

\vskip .3cm

Now an {\it arbitary} partition may be considered a dominant weight for any
$GL({\bf Z}^n)$ such that $n \geq$ the number of rows in that partition. As soon as
this condition on $n$ is met for both $\la$ and $\mu$,
$\chi(V_{\bf Z}(\mu), V_{\bf Z}(\la))$
is independent of $n$ for the following reason. Let $N > n \geq$ the number of rows
in $\la$ as well as those in $\mu$. For the moment let $\la_N$ (respectively $\la_n$)
denote the dominant weight corresponding to the partition $\la$ for the group
$GL({\bf Z}^N)$ (respectively $GL({\bf Z}^n)$) and likewise define $\mu_n, \mu_N$.
Then one has
$Ext^i_{GL({\bf Z}^n)}(V_{\bf Z}(\mu_n), V_{\bf Z}(\la_n)) \simeq
Ext^i_{GL({\bf Z}^N)}(V_{\bf Z}(\mu_N), V_{\bf Z}(\la_N))$ by, e.g., [Kulkarni2,
Proposition 1.1].
So henceforth we will deal with pairs of arbitrary partitions, always assuming that
we are working over an appropriate $GL(F)$, i.e., one for which $rank(F)$ is large
enough for the partitions in question to be valid weights.

\vskip .3cm

2.1. {\it An algorithm to compute $\chi(V_{\bf Z}(\mu), V_{\bf Z}(\la))$ for $GL(F)$.}
Let us begin by setting up an induction and by making some simple reductions.

\vskip .3cm

(1) If $\mu \not < \la$ then $\chi(V_{\bf Z}(\mu), V_{\bf Z}(\la)) = 1$ by Theorem
1.1(ii). Therefore we may assume that $\mu < \la$ and we may further induct on the
number of steps by which $\mu$ and $\la$ differ in (any linearization of) the
dominance partial order. The base case of this induction is known due to
[Kulkarni2, Theorem 2.1] and is subsumed in (2.1.3) below.

\vskip .3cm

(2) Note that $\mu < \la$ means in particular that $\la$ and $\mu$ have the same
degree, i.e., the same number of boxes. (Of course the necessity of this condition
for nontriviality of $\chi$ is clear even before by considering the action of the
center of $GL(F)$.) Let us induct on this degree as well. In degree 1 there is
only one Weyl module, namely $F$, and $\chi(F,F)$ is 1 by either part of Theorem
1.1.

\vskip .3cm

(3) If the first rows (or columns) of $\mu$ and $\la$ have the same length then by
[Kulkarni2, Proposition 1.2] we may strip them off without affecting the $Ext$ groups
and appeal to induction. So we may assume that $\mu$ and $\la$ differ in the lengths
of their first rows and also in the lengths of their first columns.

\vskip .3cm

Now let us proceed to the main inductive step. Let $\la'$ be the partition obtained
by removing the first column of $\la$ and call the length of this column $t$. By
reductions (1) and (3), the first column of $\mu$ is longer than that of $\la$.
Removing the first $t$ boxes in this column gives the skew partition $\mu/1^t$.
The partition $1^t$ is just a single column of length $t$ and the corresponding
Weyl module is $\Lambda^t(F) = $ the $t$-th exterior power of the defining
representation $F$. $\Lambda^t(F)$ is also the dual Weyl module corresponding to
 $1^t$. Consider the following equality given by the Skew Representative Theorem
from [Kulkarni1].
$$\chi(V_{\bf Z}(\mu), \Lambda^t(F) \otimes V_{\bf Z}(\la')) =
                          \chi (V_{\bf Z}(\mu/1^t), V_{\bf Z}(\la')), \leqno (2.1.1)$$
where $V_{\bf Z}(\mu/1^t)$ is the {\it skew Weyl module} corresponding to $\mu/1^t$
defined in [ABW], where it is denoted $K_{\mu/1^t}(F)$. Let us analyze both sides of
the equality (2.1.1) using the almost multiplicativity of $\chi$.

\vskip .3cm

On the right hand side, $V_{\bf Z}(\mu/1^t)$ has a Weyl filtration with factors
$V_{\bf Z}(\mu^1),\ldots,V_{\bf Z}(\mu^k)$ where $\mu^i$ are all the partitions
contained in $\mu$ such that $\mu/\mu^i$ consists of $t$ boxes no two of which
are in the same row. In other words, $\mu^i$ are all the partitions obtainable
by removing $t$ boxes from the rightmost border strip of $\mu$. This is the
characteristic-free skew Pieri rule from [AB1, Section 3], where one can also
find a more formal statement and a proof in the contravariant dual case of Schur
modules, which carries over easily to ours. By induction on the degree
$\chi(V_{\bf Z}(\mu^i), V_{\bf Z}(\la'))$ are all known.

\vskip .3cm

By the characteristic-free Pieri rule in [AB1, Section 3],
$\Lambda^t(F) \otimes V_{\bf Z}(\la')$ has a Weyl filtration with factors
$V_{\bf Z}(\la^1), V_{\bf Z}(\la^2), \ldots, V_{\bf Z}(\la^m)$, where $\la^j$
are all the partitions containing $\la'$ such that $\la^j/\la'$ consists of $t$
boxes no two of which are in the same row. Clearly these include $\la = \la^1$
(say), and $\la > \la^j$ for $j>1$. So by induction on the dominance order
$\chi(V_{\bf Z}(\mu), V_{\bf Z}(\la^j))$ are all known for $j>1$.

\vskip .3cm

{\it Case 1.} If $Hom(V_{\bf Z}(\mu), \Lambda^t(F) \otimes V_{\bf Z}(\la')) \simeq
Hom(V_{\bf Z}(\mu/1^t), V_{\bf Z}(\la'))$ is zero, i.e., if none of the $\la^j$
equals $\mu$, i.e., if none of the $\mu^i$ equals $\la'$, then $\chi$ is multiplicative
as we glue the $V_{\bf Z}(\la^j)$ together to get $\Lambda^t(F) \otimes V_{\bf Z}(\la')$
(respectively, the $V_{\bf Z}(\mu^i)$ to get $V_{\bf Z}(\mu/1^t)$) because all the $Hom$
terms in the associated long exact sequences vanish. Therefore we get the following
equations.
$$\eqalign{
\chi(V_{\bf Z}(\mu), \Lambda^t(F) \otimes V_{\bf Z}(\la')) & =
\prod  _j \chi(V_{\bf Z}(\mu),  V_{\bf Z}(\la^j)). \cr
\chi (V_{\bf Z}(\mu/1^t), V_{\bf Z}(\la')) &= \prod  _i \chi (V_{\bf Z}(\mu^i), V_{\bf Z}(\la')).\cr}
\leqno(2.1.2)$$
Now one easily finds $\chi(V_{\bf Z}(\mu), V_{\bf Z}(\la))$ from (2.1.1) and (2.1.2).

\vskip .3cm

{\it Case 2.} If Case 1 does not happen then it is clear from the Pieri rules that
$\chi$ fails to be multiplicative at exctly one stage in gluing. By the work done in
Example 2 in 1.2, the equations (2.1.2) now need to be modified by writing certain
correction factors on their right hand sides. At first glance this seems to make it
necessary to compute these integers every time we are in Case 2. The only way I
know of directly doing this computation involves finding explicit generators of certain
$Hom$ groups. This is theoretically possible, but quite hard in practice, even in the
relatively simple case treated in [Kulkarni2, Theorem 2.1], as seen in the proofs of
Lemmas A--C there.

\vskip .3cm

Fortunately we can wriggle out of this difficulty in most cases by analyzing the
combinatorics of Pieri's rule and by using conjugate symmetry of $Ext$ groups (1.6.1).
Vizualizing the dominant weights occuring in $\Lambda^t(F) \otimes V_{\bf Z}(\la')$
reveals that if one encounters Case 2 then $\mu$ must be obtainable from $\la$
by removing some of the boxes from the rightmost border strip of $\la$ and placing
them at the bottom of the first column.

\vskip .3cm

Now (1.6.1) gives
$\chi(V_{\bf Z}(\mu), V_{\bf Z}(\la)) = \chi(V_{\bf Z}(\tm), V_{\bf Z}(\tl)).$
So if we encounter Case 2, we may discard the pair $(\mu, \la)$ of dominant weights
and work with the pair $(\tl, \tm)$ instead. A little further thought shows that
the only way $(\tl, \tm)$ will also lead to Case 2 is when $\mu$ has a single box
in its last row and $\la$ is obtained by removing this box and placing it at the
end of the first row. (Recall that by the reduction (3), $\mu$ and $\la$ must differ
in their first rows as well as in their first columns.) So $\la - \mu =$ a positive
root $\alpha$. But this is precisely the situation treated by [Kulkarni2, Theorem
2.1]. This theorem gives the following calculation of all $Ext$ groups for such pairs
of Weyl modules. $Ext^1$ is cyclic of order $\langle \mu + \rho, \alpha \check{}\,\rangle +1=$
hook length of the box in the first row and first column of $\mu$ (or $\la$), and
all other $Ext$ groups vanish. So for $\la = \mu + \alpha$ one has
$$\chi(V_{\bf Z}(\mu), V_{\bf Z}(\la)) = {1 \over \langle \mu + \rho, \alpha \check{}\,\rangle +1}.
\leqno(2.1.3)$$
This ends the recursion and completes the description of the algorithm.

\vskip .3cm

{\it Remarks.} (1) For future use note the following extra flexibility that may be built
into the above algorithm. $t$ may be chosen to be the length of any column in $\la$.
Then $\la'$ would be the partition obtained by deleting the rightmost box in each
of the top $t$ rows of $\la$. The induction goes through just as before with the
following differences. For one, different sets of $\mu^i$ and $\la^j$ are involved in
the computation. Moreover the characteriztion of when one encounters Case 2 is no longer
as clean. Nor is the escape by use of conjugate symmetry guaranteed if one is required
to use a given $t$. But of course one may then use a different $t$ (e.g., just follow
the algorithm above). The point is that sometimes choosing a different $t$ will make
the calculation simpler. For example choosing $t=1$ (of course the last column of $\la$
must be of length 1 for that) may considerably simplify the combinatorics for large
partitions. Note that [Kulkarni2, Theorem 2.1] was proved by taking $t=1$. In the proof
of that theorem recursion could be controlled sufficiently to give all $Ext$ groups, not
just $\chi$.

\vskip .3cm

(2) (Relating the algorithm to symmetries of $Ext$ groups between Weyl modules.)
By (1.6.1), (1.6.2) and using the terminology there we have
$$\chi (V_{\bf Z}(\mu^c), V_{\bf Z}(\la^c)) = \chi(V_{\bf Z}(\mu), V_{\bf Z}(\la)) =
\chi (V_{\bf Z}(\tl), V_{\bf Z}(\tm)).$$
Under the isomorphism (1.6.2), using the algorithm for the pair of partitions
$(\mu,\la)$ by splitting off the first column of $\la$ is seen to be tantamount
to using the algorithm for the pair $(\mu^c,\la^c)$ by splitting off the last column
of $\la^c$. Similarly one can devise a version of the algorithm that is consistent
with the second equality. In place of (2.1.1) this version will rely on a procedure
that splits off the top row of $\mu$ rather than the first column of $\la$. Later
both of these symmetries will be very useful for substantial reductions while
proving a general formula for $\chi(V_{\bf Z}(\mu), V_{\bf Z}(\la))$.

\vskip .3cm

2.2. {\it Examples.} Before stating and proving the formula for $\chi$ for a general
pair of partitions, it will be instructive to see the above algorithm in action in
some simple examples.

\vskip .3cm

{\it Example 1.} ($GL_2.)$ Let us compute $\chi(V_{\bf Z}(\mu), V_{\bf Z}(\la))$,
where $\mu$ and $\la$ are partitions with at most two parts. By the initial reductions
the general case immediately reduces to the following: $\mu= (a,b)$, $\la=(a+b,0)$,
where $a \geq b >0$. For convenience let us denote $\chi(V_{\bf Z}(\mu), V_{\bf Z}(\la))$
by $\chi[a,b]$. Clearly the only choice for $t$ here is $t=1$. Unless $b=1$, the algorithm
does not lead to Case 2. Using (2.1.1), (2.1.2) and reduction (3) one gets the following
equations.
$$\eqalign{
\chi[a,b] \; \chi[a-1, b-1] &= \chi[a,b-1] \; \chi[a-1, b] \quad {\rm if } \; a>b.\cr
\chi[a,a] \, \chi[a-1, a-1] &= \chi[a,a-1].}$$
If $b=1$ one has $\chi[a,1] = {1 \over a+1}$. Now one easily derives that
$$\chi[a,b] = {b \over a+1}.$$

{\it Example 2.} ($\chi$ for two hook partitions.) Let $\mu = (a,1^b), \la = (a+s,1^{b-s})$,
i.e., both diagrams have the shape of a hook. For convenience let us denote
$\chi(V_{\bf Z}(\mu), V_{\bf Z}(\la))$ by $\chi[a,b,s]$. As in the algorithm, let
us choose $t=b+1-s$, i.e. the length of the first column of $\la$ (The reason for
this choice will be clear soon.) If $s \neq 1$, the algorithm does not lead to Case 2.
In this case one gets the following equation using (2.1.1) and (2.1.2).
$$\chi[a,b,s] \; \chi[a,b,s-1] = \chi[a,s-1,s-1] \; \chi[a-1, s,s].$$
If $s=1$ one has $\chi[a,b,1]= {1 \over a+b}$. Now one easily derives that
$$\chi[a,b,s] = \left({a+b \over s}\right)^{(-1)^s}.$$
Notice that for this calculation, we only needed to deal with pairs of hooks
i.e., no other partitions appeared in recursion. A little thought reveals that
if we had chosen $t=1$, this would not have been the case and the computation
would have been much more unwieldy.

\vskip .3cm

2.3. {\it A formula for $\chi(V_{\bf Z}(\mu), V_{\bf Z}(\la))$.} It turns out that
there is a fairly simple formula for $\chi(V_{\bf Z}(\mu), V_{\bf Z}(\la))$ in
terms of the geometry of the diagrams of the involved partitions. Once the formula
is guessed (after working out several classes of examples like the ones shown above),
it is not so hard to prove by judicious use of the algorithm and careful bookkeeping.
Before stating the final answer in the next theorem, let us set up some terminology.
A skew partition is {\it connected} if a rook can go from any box in it to any
other by ordinary chess moves. A {\it skew hook} (also called ribbon or a border
strip) is a skew partition not containing a 2 by 2 square. We may call a skew
partition failing this condition {\it overconnected}. Thus a skew partition will
fail to be a connected skew hook by being disconnected or overconnected or both.
We could call a connected skew hook a {\it snake}.
The ``right endpoint" (respectively, ``left endpoint") of a skew partition will mean
the rightmost box in its top row (respctively, the leftmost box in its bottom row).
This terminology will be used mainly when the skew partition is a skew hook.

\vskip .3cm

{\it Theorem.} Let $\mu$ and $\la$ be arbitrary partitions. Considering
$\mu$ and $\la$ as sets of appropriately situated boxes in a plane, let $\nu$ be
the partition $\la \cap \mu$. We will say that $\mu$ and $\la$ differ by
connected skew hooks if the skew partitions $\mu/\nu$ and $\la/\nu$ are both
connected skew hooks. Then one has the following.
(1) $\chi(V_{\bf Z}(\mu), V_{\bf Z}(\la))$ is 1 unless $\mu < \la$ and moreover
$\mu$ and $\la$ differ by connected skew hooks.
(2) If $\mu$ and $\la$ do differ by connected skew hooks and $\mu < \la$ as well,
then
$$\chi(V_{\bf Z}(\mu), V_{\bf Z}(\la)) =  \left({\ell \over d}\right) ^ {(-1) ^ r},$$
where the symbols have meanings as follows. $\ell$ is the the equal number of boxes in
skew hooks $\mu/\nu$ and $\la/\nu$. We will call this the {\it length} of the skew hook,
not to be confused with the number of rows in the skew hook. $d$ is the total
number of right and up moves that each box in $\mu/\nu$ has to make (i.e., the distance
this snake has to cover by sliding along the border of the diagram) in order to transform
$\mu$ into $\la$. $r$ is the sum of the number of nonempty rows in $\mu/\nu$ and that
in $\la/\nu$. By drawing pictures, the numbers $\ell$ and $d$ can also be seen to be the
following hook lengths. $\ell$ is the hook length of the box in $\la$ that is in the same row
as the right endpoint of $\la/\nu$ and the same column as the left endpoint of $\la/\nu$. The
previous sentence stays valid after replacing each $\la$ by $\mu$. $d$ is the hook length
of the box in $\mu$ that is in the same row as the left endpoint of $\la/\nu$ and the same
column as the left endpoint of $\mu/\nu$. $d$ is also the hook length of the box in $\la$
that is in the same row as the right endpoint of $\la/\nu$ and the same column as the right
endpoint of $\mu/\nu$.

\vskip .3cm

{\it Proof}.
Recall the following from 2.1. By Theorem 1.1(ii) we may assume that $\mu < \la$. We may
also assume that $\la$ has a strictly longer first row and strictly shorter first column
than $\mu$. So the entire last column of $\la$ is missing in $\mu$ and the entire last
row of $\mu$ is missing in $\la$. We will induct on the equal number of boxes in $\mu$
and $\la$ and on the number of steps by which $\la$ and $\mu$ differ under (a linearization
of) the dominance partial order. Note that the base case of each induction is already
known. By (2.1.3) we may additionally assume that $\la/\nu$ (and $\mu/\nu$) contains
more than one box. We will follow the algorithm as extended in the first remark following
it. More precisely we will always use $t =$ the length of the {\it last} column of
$\la$. It will be clear that the difficulty with infinite $Hom$ encountered in Case 2
in 2.1 never arises (partly due to the use of (2.1.3) at the outset).

\vskip .3cm

The proof will consist of a series of reductions and one crucial calculation
requiring treatment of several cases. Throughout it will be easier to
follow the combinatorial arguments by visualizing the diagrams of $\la$ and $\mu$
with $\la/\nu$ and $\mu/\nu$ (the set differences between $\la$ and $\mu$) ``colored"
differently from $\nu$.

\vskip .3cm

{\it Step 1.} (Getting rid of almost all overconnectedness.)
Suppose that $\la/\nu$ contains a 2 by 2 square (i.e., is overconnected) and that
$\la$ has at least one column strictly to the right of this 2 by 2 square. When one
runs the main step of the algorithm with $t$ = the length of such a column (say the last
one), for each resulting pair $(\mu,\la^j)$ the skew partition $\la^j/(\la^j \cap \mu)$
will stay overconnected. This is simply because $\la^j$ is obtained by moving down some
boxes to the right of the ``excess 2 by 2 square," which will clearly keep this square
intact. Similarly for each pair $(\mu^i,\la')$ the skew partition $\la'/(\la' \cap \mu^i)$
will also stay overconnected. Now equations (2.1.1), (2.1.2) and induction give the desired
triviality of $\chi(V_{\bf Z}(\mu), V_{\bf Z}(\la))$. By conjugate symmetry (1.6.1), one
also gets the triviality of $\chi(V_{\bf Z}(\mu), V_{\bf Z}(\la))$ when $\mu/\nu$ contains
a 2 by 2 square and $\mu$ has at least one row strictly below this 2 by 2 square.

\vskip .3cm

Enclose $\mu$ and $\la$ in a rectangle and recall the symmetry of $Ext$ under complements
(1.6.2). Since any 2 by 2 square in $\la/\nu$ is present is $\la$ and absent in $\mu$,
the same square will be present in $\mu^c$ and absent in $\la^c$. If such an excess 2 by
2 square occurs in $\la$ and is situated below the top row of $\la$, the corresponding
excess square in $\mu^c$ must be situated above the bottom row of $\mu^c$. Now applying
the last sentence of the previous paragraph to the pair $(\mu^c, \la^c)$ and using (1.6.2)
we get that $\chi(V_{\bf Z}(\mu), V_{\bf Z}(\la)) = 1$. Conjugate symmetry (1.6.1) gives
the same conclusion if $\mu/\nu$ contains a 2 by 2 square that is situated to the right
of the first column of $\mu$.

\vskip .3cm

By the previous two paragraphs we may assume the following. $\la/\nu$ contains at most
one 2 by 2 square and such a square must occur in the last two columns and the first two
rows of $\la$. Similarly $\mu/\nu$ contains at most one 2 by 2 square and such a square
must occur in the first two columns and the last two rows of $\mu$. Henceforth these
assumptions will be in force throughout the proof.

\vskip .3cm

{\it Step 2.} (The crucial calculation.) Suppose the last column of $\la$ (which is absent
from $\mu$ by our reduction) consists of exactly one box. We will run the algorithm with
$t = 1$ and show by induction that the claimed formula holds. By Step 1, $\la/\nu$ is a
(possibly disconnected) skew hook. We will treat three cases.

\vskip .3cm

{\it Case A.} Suppose $\la/\nu$ and $\mu/\nu$ are both connected skew hooks. We will
analyze both equations in (2.1.2). In the first equation it is easy to see by induction
that there are at most three partitions $\la^j$ other than $\la$ that lead to nontrivial
$\chi(V_{\bf Z}(\mu), V_{\bf Z}(\la^j))$. These are described in cases (A.1.i) through
(A.1.iii) below. For ease of description let the right endpoint of $\mu/\nu$ be in $u$-th
row and $v$-th column of $\mu$. Let the left endpoint of $\la/\nu$ is in $x$-th row and
$y$-th column of $\la$.

\vskip .3cm

(A.1.i) When $\la^j$ is obtained by adding a box to $\la'$ immediately below the left
endpoint of $\la/\nu$, i.e., in $(x+1)$-th row and $y$-th column. This gives a valid
partition precisely when the $x$-th and $(x+1)$-th rows of $\nu$ are of equal length $y-1$.
In that case $\la^j \cap \mu =\nu$ and moreover $\la^j$ and $\mu$ still differ by connected
skew hooks of length $\ell$. But the distance needed to slide one of these skew hooks into
another is one less than before and $\la^j/\nu$ has one more row than $\la/\nu$. By
induction
$$\chi(V_{\bf Z}(\mu), V_{\bf Z}(\la^j)) = \left({\ell \over {d-1}} \right)^{(-1)^{r+1}}.$$

(A.1.ii) When $\la^j$ is obtained by adding a box to $\la'$ in the position of the left
endpoint of $\mu/\nu$. This gives a valid partition precisely when $\mu$ has exactly
one row more than $\la$. In a manner similar to (A.1.i) it is easy to see that in this
case we have
$$\chi(V_{\bf Z}(\mu), V_{\bf Z}(\la^j)) = \left({{\ell - 1} \over {d-1}} \right)^{(-1)^{r}}.$$

(A.1.iii) When $\la^j$ is obtained by adding a box to $\la'$ in the position of the right
endpoint of $\mu/\nu$, i.e., in $u$-th row and $v$-th column. This procedure gives a valid
partition precisely when in $\mu/\nu$ this right endpoint is the only box in its row,
i.e., when the $u$-th row of $\nu$ has length $v-1$. In a manner similar to (A.1.i) and
(A.1.ii) it is easy to see that in this case we have
$$\chi(V_{\bf Z}(\mu), V_{\bf Z}(\la^j)) = \left({{\ell-1} \over {d}} \right)^{(-1)^{r-1}}.$$

The analysis of the second equation is very similar. Here we need to look for the $\mu^i$
that will give nontrivial $\chi(V_{\bf Z}(\mu^i), V_{\bf Z}(\la'))$. There are at most three
possibilities, which are listed below.

\vskip .3cm

(A.2.i) When $\mu^i$ is obtained from $\mu$ by removing the box immediately to the left of
the left endpoint of $\la/\nu$, i.e., the one in $x$-th row and $(y-1)$-th column. This gives
a valid partition precisely when $x$-th and $(x+1)$-th rows of $\nu$ are of unequal legnths,
i.e., when (A.1.i) is not possible. In that case $\la'$ and $\mu^i$ still differ by connected
skew hooks of length $\ell$; the distance needed to slide one of these skew hooks into
another is one less than before; and the number of rows in each skew hook stays unchanged.
By induction
$$\chi(V_{\bf Z}(\mu^i), V_{\bf Z}(\la')) = \left({\ell \over {d-1}} \right)^{(-1)^{r}}.$$

(A.2.ii) When $\mu^i$ is obtained from $\mu$ by removing the box at the left endpoint of
$\mu/\nu$. It is easy to see that this gives a valid partition precisely when (A.1.ii) is not
possible. Moreover, in that case $\chi(V_{\bf Z}(\mu^i), V_{\bf Z}(\la'))$ is the reciprocal
of the number obtained in (A.1.ii).

\vskip .3cm

(A.2.iii) When $\mu^i$ is obtained from $\mu$ by removing the box at the right endpoint
of $\mu/\nu$. Again this gives a valid partition precisely when (A.1.iii) is not possible. In that
case $\chi(V_{\bf Z}(\mu^i), V_{\bf Z}(\la'))$ is the reciprocal of the number obtained in (A.1.iii).

\vskip .3cm

Note that a fourth $\mu^i$ seems possible at first glance. Namely one could try to take off a box
from $\mu$ so as to somehow augment the skew hook $\la/\nu$ at its right endpoint (similar to
the way it was augmented on the left in (A.2.i)). But this is clearly seen to be impossible. We
will encounter such a possibility later in Case C.

\vskip .3cm

So we may as well suppose that none of the cases (A.2.i) through (A.2.iii) occurs and thus all
of the cases (A.1.i) though (A.1.iii) do. The desired claim is now immediate from (2.1.1) and
the following trivial calculation.
$$\left({\ell \over {d}} \right)
\left({{d-1} \over \ell} \right)
\left({{\ell - 1} \over {d-1}} \right)
\left({d \over {\ell-1}} \right) = 1.$$

{\it Case B.} Suppose $\la/\nu$ is a connected skew hook but $\mu/\nu$ is not a connected skew
hook. Proceeding in a manner similar to Case A, it is easy to see that there are only two ways
in which one of the equations in (2.1.2) could yield a lower term with nontrivial $\chi$. (If
every lower $\chi$ in both equations is trivial, the desired result, namely triviality of
$\chi(V_{\bf Z}(\mu), V_{\bf Z}(\la))$, is immediate.)

\vskip .3cm

(B.1) If $\mu/\nu$ contains a 2 by 2 square (by Step 1 necessarily only one, occuring in
the last two rows and the first two columns of $\mu$), then in each equation there is at
most one nontrivial term possible. In the first equation this occurs when $\la^j$ is
obtained by adding a box to $\la'$ in the position of the top left box of the 2 by 2
square. In the second equation this occurs when $\mu^i$ is obtained from $\mu$ by removing
the bottom right box of this 2 by 2 square. These two possibilities are easily seen to
give the same value of $\chi$. So they cancel each other in (2.1.1) and give the desired
triviality of $\chi(V_{\bf Z}(\mu), V_{\bf Z}(\la))$.

\vskip .3cm

(B.2) If $\mu/\nu$ does not contain a 2 by 2 square, it must be a disconnected skew hook.
It is easy to see that the only way to get a nontrivial term in the recursion is when
$\mu/\nu$ has two connected components and one of the components is a single box. In this
case each equation has exactly one nontrivial lower term and again these two terms cancel
each other. (In the first equation get $\la^j$ by adding a box to $\la'$ in place of the
isolated box in $\mu/\nu$. In the second equation get $\mu^i$ by removing the same isolated
box from $\mu$.)

\vskip .3cm

{\it Case C.} Suppose $\la/\nu$ is a disconnected skew hook. Proceeding as before, we will
analyze the ways in which one of the equations in (2.1.2) could yield a nontrivial lower $\chi$.

\vskip .3cm

(C.1) If $\la'/\nu$ is a disconnected skew hook as well, then it is easy to see that the only
way one could get nontrivial lower terms in either equation is when the all of the following
three conditions hold. $\la'/\nu$ must consist of exactly two connected components. Moreover
the gap between the components must be exactly one box, i.e., adding just one box to $\la'/\nu$,
say in $a$-th row and $b$-th column of $\la'$, should make it a connected skew hook. And
$\mu/\nu$ must be a connected skew hook. In that case there is exactly one nontrivial term
in each equation and once again these cancel each other. ($\la^j$ is obtained by adding a
box to $\la'$ in $a$-th row and $b$-th column. In the second equation obtain $\mu^i$ from
$\mu$ by deleting the box in $(a-1)$-th row and $(b-1)$-th column.)

\vskip .3cm

(C.2) So suppose now that $\la'/\nu$ is a connected skew hook (i.e., $\la/\nu$ has exactly two
connected components and one of the components consists of the last box in the first row of
$\la$.)  Now if $\mu/\nu$ fails to be a connected skew hook, it is easy to see that the entire
analysis in Case B carries over. (After changing $\la$ to $\la'$ in the opening sentence of
Case B, the rest applies verbatim.)

\vskip .3cm

(C.3) So we may suppose that $\mu/\nu$ and $\la'/\nu$ are both connected skew hooks.  The
analysis here is entirely parallel to that in Case A with the following crucial difference.
There is a fourth nontrivial term possible in each equation in (2.1.2). In the first
equation one could get a $\la^j$ by adding a box to $\la'$ immediately to the right of the
right endpoint of $\la'/\nu$. (In Case A this would just give $\la^j = \la$, but here we get
a lower term.) In the second equation one could get another $\mu^i$ from $\mu$ by deleting
the box immediately above the right endpoint of $\la'/\nu$. (This is the putative fourth
possibility discussed above immediately after (A.2.iii), which could not occur there.) Further,
just as for the three pairs of possibilites in Case A, exactly one of this new pair of
possibilities will actually take place. It is easy to see that the same calculation as in
Case A gives the desired triviality of $\chi(V_{\bf Z}(\mu), V_{\bf Z}(\la))$.

\vskip .3cm

Clearly cases A, B and C exhaust all possibilities when the last column of $\la$ has
length 1 (in presence of the reductions that were made previously and which will be in
force throughout the proof). This finishes Step 2.

\vskip .3cm

{\it Step 3.} (Further reductions.)
If the last row of $\mu$ (which is absent from $\la$ by earlier reduction) contains
exactly one box we will be done by Step 2 and conjugate symmetry (1.6.1). So henceforth
we will assume that the last column of $\la$ as well as the last row of $\mu$ (each of which
is entirely missing from the other partition) have lengths greater than 1.

\vskip .3cm

Now suppose that the first column of $\mu$ has exactly one more box than the
first column of $\la$. Then one reduces to the case in Step 2 as follows. In view
of the discussion near the beginning of Section 2 we may take $n$ (the rank of the
defining representation of $GL(F)$) to be the number of rows in $\mu$. Now by (1.6.2)
one may replace the pair of partitions $(\mu, \la)$ by the pair $(\mu^c, \la^c)$,
which is covered by Step 2. So henceforth we will assume that the first column of
$\mu$ has at least two boxes more than the first column of $\la$, i.e., that the
entire last two rows of $\mu$ are missing in $\la$. By conjugate symmetry (1.6.1)
we may assume that the entire last two columns of $\la$ are missing in $\mu$.

\vskip .3cm

{\it Step 4.} (The last case.)
Combining Step 1 and Step 3 we may and will assume the following. Each of
$\la/\nu$ as well as $\mu/\nu$ contains exactly one 2 by 2 square and these
squares are situated as described at the end of Step 1. In particular note that
the last column of $\la$ contains exactly two boxes, as does the last row of
$\mu$. We will run the algorithm in 2.1 with $t = 2$ and prove that
$\chi(V_{\bf Z}(\mu), V_{\bf Z}(\la))$ is 1.

\vskip .3cm

Consider all pairs of partitions $(\mu, \la^j)$ (except $(\mu,\la)$) that occur
in the first equation in (2.1.2). For such a pair to lead to nontrivial $\chi$,
one of the two removed boxes in the last column of $\la$ must be added to the
bottom of the first column of $\la'$. This is the only way to remove the
overconnectedness in $\la/\nu$ and $\mu/\nu$. So essentially one has to deal with
placement of just one box. Similarly for a pair $(\mu^i, \la')$ to lead to nontrivial
$\chi$ in the second equation in (2.1.2), the second box in the last row of $\mu$
must be missing in $\mu^i$. Thus essentially one has to deal with removal of only
one box in $\mu$. Thus we have a situation very similar to the one handled in Step
2 and one can imitate the cases there. Instead let us take a shortcut.

\vskip .3cm

For future need let us separately treat the case when each of $\la/\nu$ and $\mu/\nu$
contains exactly four boxes, i.e., is a 2 by 2 square. Using preceding discussion,
one easily verifies the desired result by an explicit calculation that is similar to
but simpler than Case A in Step 2. So from now on suppose that each of $\la/\nu$ and
$\mu/\nu$ contains more than four boxes.

\vskip .3cm

Let $\xi$ be the partition obtained from $\la$ by removing the two boxes in the
last column and adding one box at the bottom of the first column. One has
$$\chi(V_{\bf Z}(\mu), V_{\bf Z}(\la') \otimes \Lambda^2(F)) =
  \chi(V_{\bf Z}(\mu), V_{\bf Z}(\la)) \; \chi(V_{\bf Z}(\mu), V_{\bf Z}(\xi) \otimes F)$$
for the following reason. Recasting earlier discussion, any lower $\la^j$ that leads to
nontrivial $\chi(V_{\bf Z}(\mu), V_{\bf Z}(\la^j)$ must contain $\xi$ as a subset. Thus,
apart from $\la^1 = \la$, it suffices to consider only those $\la^j$ that are obtained by
adding one box to $\xi$. This is clearly a subset of the set of partitions $\xi^k$ such
that $V_{\bf Z}(\xi^k)$ is a filtration factor in a Weyl filtration of $V_{\bf Z}(\xi) \otimes F$.
In fact the two sets will be the same unless the first two columns of $\la$ have equal
lengths, in which case there will be exactly one extra $\xi^k$ obtained by adding a box
in the last row of $\xi$. (This partition will not occur among the $\la^j$ since it amounts
to adding two boxes to $\la'$ in the same row.) However in this case $\mu/(\mu \cap \xi^k)$
will be disconnected since one of the connected components will be the two boxes in the last
row and there must be at least three boxes in $\mu/(\mu \cap \xi^k)$ thanks to the case
treated separately in the previous paragraph. Thus for the extra $\xi^k$, one has
$\chi(V_{\bf Z}(\mu), V_{\bf Z}(\xi^k) = 1$. This proves the claimed equality.

\vskip .3cm

Now using (2.1.1) one has
$$\eqalign{
\chi(V_{\bf Z}(\mu), V_{\bf Z}(\la') \otimes \Lambda^2(F)) & =
                                               \chi(V_{\bf Z}(\mu/1^2), V_{\bf Z}(\la')).\cr
\chi(V_{\bf Z}(\mu), V_{\bf Z}(\xi) \otimes F) & =
                                               \chi(V_{\bf Z}(\mu/1), V_{\bf Z}(\xi)).\cr}$$
We will be done once we show that the right hand sides of the preceding two equations
are the equal. To see this let $\pi$ be the partition obtained by deleting the
second box from the last row of $\mu$. Recasting earlier discussion, any $\mu^i$ that
leads to nontrivial $\chi(V_{\bf Z}(\mu^i), V_{\bf Z}(\la')$ must be obtainable by removing
one box from $\pi$. Thus while calculating $\chi(V_{\bf Z}(\mu/1^2), V_{\bf Z}(\la'))$ it
suffices to consider only those $\mu^i$ that are obtained by deleting one box from $\pi$.
We will set a up an {\it almost} pairing between the set $S$ of such partitions and the set
$T$ of all partitions $\bar\mu$ such that $V_{\bf Z}(\bar \mu)$ is a filtration factor in a
Weyl filtration of $V_{\bf Z}(\mu/1)$. Suppose $\mu^i \in S$ is obtained by removing a
certain box from $\pi$. Define a map $f$ from $S$ to $T$ by sending $\mu^i$ to the
partition obtained by removing the same box from $\mu$. Before examining the failure
of $f$ to be a pairing, let us note that
$$\chi(V_{\bf Z}(\mu^i), V_{\bf Z}(\la')) = \chi(V_{\bf Z}(f(\mu^i)), V_{\bf Z}(\xi)).$$
To see this one easily verifies pictorially that $\la'/(\la' \cap \mu^i) = \xi/(\xi \cap f(\mu^i))$
and that the only difference between $\mu^i /(\la' \cap \mu^i)$ and $f(\mu^i) /(\xi \cap f(\mu^i))$
is in the arrangement of the three boxes at the left endpoint, where the former has the same
shape as the partition (2,1) and the latter has the same shape as the skew partition (2,2)/(1).
Clearly this difference is irrelevant for the formula asserted in Theorem 2.3. Further, the set
differences between each pair of partitions have the same relative position. The stated equality
follows by induction on the degree.

\vskip .3cm

Now let us deal with the failure of $f$ to be a pairing. This can happen in two ways.
(1) $f$ will not yield the partition $\bar\mu$ in $T$ that is obtained by removing the second box
in the last row of $\mu$. However for this partition $\chi(V_{\bf Z}(\bar\mu), V_{\bf Z}(\xi))$
is trivial since $\bar\mu/(\bar\mu\cap \xi)$ is disconnected in the last two rows. So we may
ignore this failure. (2) $f(\mu^i)$ will be undefined if the last but one row of $\mu$ consists
of exactly two boxes and $\mu^i$ is obtained by removing the second of these from $\pi$.
(Removing this box from $\mu$ will not give a partition.) However in this case
$\chi(V_{\bf Z}(\mu^i), V_{\bf Z}(\la'))$ is trivial since $\mu^i/(\mu^i \cap \la')$, which contains
at least three boxes thanks to the special case trated above, is disconnected. So we may
ignore this failure as well. This finishes the proof of Theorem 2.3.

\vskip .3cm

{\it Remarks.} 1) Certain products of the numbers on the right hand side of the formula
in Theorem 2.3 have already appeared in the literature on Jantzen's sum formula. [JM]
shows that the determinant of the Gram matrix for a Specht module for the symmetric
group is such a product. (I am grateful to Arun Ram for pointing out this fact and
for supplying the reference.)

\vskip .3cm

2) Brenti's recent work on certain parabolic Kazhdan-Lusztig and $R$-polynomials for the
symmetric group also involves connected skew hooks, see [Brenti]. Theorem 2.3 and Brenti's
work together suggest that at least for type A, $Ext$ groups between Weyl modules should be
somehow related to Kazhdan-Lusztig combinatorics. In the BGG category a connection between
ordinary $R$-polynomials for the Weyl group and $Ext$ groups between Verma modules was
suggested by Gabber and Joseph, but their guess was found to be false by Boe. Nonetheless
in light of new evidence it seems that there may well be a relationship between appropriate
$R$-polynomials and $Ext$ groups between Weyl modules. It would be very interesting to
find a precise connection.

\vskip 0.8cm

\noindent {\it References}

\vskip .3cm

\parindent = 0pt

[AB1] K. Akin and D. A. Buchsbaum, Characteristic-free representation
theory of the general linear group,
{\it Adv. in Math.} {\bf 58} (1985), 149--200.

\vskip .2cm

[AB2] K. Akin and D. A. Buchsbaum, Characteristic-free representation
theory of the general linear group II: Homological considerations,
{\it Adv. in Math.} {\bf  72} (1988), 171--210.

\vskip .2cm

[ABW] K. Akin, D. A. Buchsbaum, and J. Weyman, Schur functors and Schur complexes,
{\it Adv. in Math.} {\bf 44} (1982), 207--278.

\vskip .2cm

[AR] A. M. Adamovich, G. L. Rybnikov, Tilting modules for classical groups and Howe
duality in positive characteristic, {\it Transform. Groups} {\bf 1} (1996), 1--34.

\vskip .2cm

[Brenti] F. Brenti, Kazhdan-Lusztig and $R$-polynomials, Young's lattice, and Dyck
partitions, {\it Pacific J. Math.} {\bf 207} (2002), 257--286.

\vskip .2cm

[CPSvdK] E. Cline, B. Parshall, L. Scott, and W. van der Kallen,
Rational and generic cohomology, {\it Invent. Math.} {\bf 39} (1977), 143--163.

\vskip .2cm

[JM]  G. D. James, G. E. Murphy, The determinant of the Gram matrix for a Specht module,
{\it J. Algebra} {\bf 59} (1979), 222--235.

\vskip .2cm

[Jantzen] J. C. Jantzen, Representations of algebraic groups, second edition,
Mathematical Surveys and Monographs {\bf 107}, Amer. Math. Soc., 2003.

\vskip .2cm

[Kulkarni1] U. Kulkarni, Skew Weyl modules for ${\rm GL}_n$ and degree reduction for Schur algebras,
{\it J. Algebra} {\bf 224} (2000), 248--262.

\vskip .2cm

[Kulkarni2] U. Kulkarni, On the Ext groups between Weyl modules for ${\rm GL}_n$,
{\it preprint}, arXiv:math.RT/0505370.

\vskip .2cm

[Kulkarni3] U. Kulkarni, A proof of Andersen's tilting sum formula, {\it in preparation}.

\vskip .2cm

[McNinch] G. McNinch, Filtrations and positive characteristic Howe duality,
{\it Math Z.} {\bf 235} (2000), 651--685.

\end